\documentclass[12pt]{article}
\usepackage{amsthm,amsfonts,amssymb,amscd}

\begin{document}

\begin{center}
\large \bf Canonical and log canonical thresholds \\ of multiple
projective spaces
\end{center}\vspace{0.5cm}

\centerline{A.V.Pukhlikov}\vspace{0.5cm}

\parshape=1
3cm 10cm \noindent {\small \quad\quad\quad \quad\quad\quad\quad
\quad\quad\quad {\bf }\newline In this paper we show that the
global (log) canonical threshold of $d$-sheeted covers of the
$M$-dimensional projective space of index 1, where $d\geqslant 4$,
is equal to one for almost all families (except for a finite set).
The varieties are assumed to have at most quadratic singularities,
the rank of which is bounded from below, and to satisfy the
regularity conditions. This implies birational rigidity of new
large classes of Fano-Mori fibre spaces over a base, the dimension
of which is bounded from above by a constant that depends
(quadratically) on the dimension of the fibre only.\vspace{0.1cm}

Bibliography: 12 items.} \vspace{1cm}

\vspace{1cm}

\noindent 14E05, 14E07\vspace{1cm}

\section*{Introduction}

{\bf 0.1. Statement of the main results.} In \cite{Pukh19a}
general $d$-sheeted covers of the complex projective space
${\mathbb P}={\mathbb P}^M$ which are Fano varieties of index 1
with at most quadratic singularities, the rank of which is bounded
from below, were shown to be birationally superrigid. In this
paper we will prove that for almost all values of the discrete
parameters defining these varieties a general multiple projective
space of index 1 satisfies a much stronger property: its global
canonical (and the more so, log canonical) threshold is equal to
1. Now \cite{Pukh15a} immediately implies the birational rigidity
type results for fibre spaces, the fibres of which are multiple
projective spaces, and new classes of Fano direct products
\cite{Pukh05}. Let us give precise statements.\vspace{0.1cm}

Fix a pair of positive integers $(d,l)\in {\mathbb Z}^{\times
2}_+$ in the set described by the following table:\vspace{0.5cm}

\begin{center}
\begin{tabular}{|c|c|}
\hline
$d$  & $l$ \\
\hline 4 & $\geqslant 21$ \\
\hline 5 & $\geqslant 5$ \\
\hline 6 & $\geqslant 6$ \\
\hline 7,8 & $\geqslant 4$ \\
\hline 9,10 & $\geqslant 3$ \\
\hline $\geqslant 11$ & $\geqslant 2$ \\
\hline
\end{tabular}
\end{center}
\vspace{0.5cm}

Set $M=(d-1)l$. The symbol ${\mathbb P}$ stands for the complex
projective space ${\mathbb P}^M$. Consider the weighted projective
space
$$
\overline{{\mathbb P}}= {\mathbb
P}(\underbrace{1,\dots,1}_{M+1},l)={\mathbb P}(1^{M+1},l)
$$
with homogeneous coordinates $x_0,\dots,x_M$, $\xi$, where $x_i$
are of weight 1 and $\xi$ is of weight $l$, and a
quasi-homogeneous polynomial
$$
F(x_*,\xi)=\xi^d+A_1(x_*)\xi^{d-1}+\dots+A_d(x_*)
$$
of degree $dl$ (that is, $A_i(x_0,\dots,x_M)$ is a homogeneous
polynomial of degree $il$ для $i=1,\dots,d$). The space
$$
{\cal F}=\prod^d_{i=1}H^0({\mathbb P},{\cal O}_{\mathbb P}(il))
$$
parameterizes all such polynomials. If the hypersurface
$$
V=\{F=0\}\subset\overline{\mathbb P}
$$
has at most quadratic singularities of rank $\geqslant 7$ (and we
will consider hypersurfaces with stronger restrictions for the
rank), then $V$ is a factorial variety with terminal
singularities, see \cite{Pukh19a}, so that
$$
\mathop{\rm Pic}V={\mathbb Z}H,\quad K_V=-H,
$$
where $H$ is the class of a ``hyperplane section'', that is, of
the divisor $V\cap \{\lambda=0\}$, where $\lambda(x_0,\dots,x_M)$
is an arbitrary linear form. Below for all the values of $d,l$
under consideration we will define explicitly a positive
integral-valued function $\varepsilon(d,l)$, which behaves as
$\frac12 M^2$ as the dimension $M$ grows.\vspace{0.1cm}

As in \cite{Pukh19a}, we identify the polynomial $F\in{\cal F}$
and the corresponding hypersurface $\{F=0\}$, which makes it
possible to write $V\in{\cal F}$. The following theorem is the
main result of the present paper.\vspace{0.1cm}

{\bf Theorem 0.1.} {\it There is a Zariski open subset ${\cal
F}_{\rm reg}\subset{\cal F}$ such that:\vspace{0.1cm}

{\rm (i)} every hypersurface $V\in{\cal F}_{\rm reg}$ has at most
quadratic singularities of rank $\geqslant 8$ and for that reason
is a factorial Fano variety of index 1 with terminal
singularities,\vspace{0.1cm}

{\rm (ii)} the inequality
$$
\mathop{\rm codim}(({\cal F}\backslash{\cal F}_{\rm reg})
\subset{\cal F})\geqslant \varepsilon(d,l)
$$
holds,\vspace{0.1cm}

{\rm (iii)} for every variety $V\in{\cal F}_{\rm reg}$ and every
divisor $D\sim nH$ the pair $(V,\frac{1}{n}D)$ is
canonical.}\vspace{0.1cm}

Now \cite[Theorem 1.1]{Pukh15a} makes it possible to describe the
birational geometry of Fano-Mori fibre spaces, the fibres of which
are multiple projective spaces of index 1.\vspace{0.1cm}

Let $\overline{\eta}\colon {\mathbb X}\to S$ be a locally trivial
fibre space, the base of which is a non-singular projective
rationally connected variety $S$ of dimension
$$
\mathop{\rm dim\,} S< \varepsilon(d,l),
$$
and the fibre is the weighted projective space $\overline{{\mathbb
P}}$. Consider an irreducible hypersurface $W\subset {\mathbb X}$,
such that for every point $s\in S$ the intersection
$$
\overline{\eta}^{-1}(s)\cap W\in {\cal F}
$$
is a multiple projective space of the type described above. The
claim (ii) of Theorem 0.1 implies that we may assume that
$$
W_s=\overline{\eta}^{-1}(s)\cap W\in {\cal F}_{\rm reg}
$$
for every point of the base $s\in S$, if the linear system $|W|$
is sufficiently mobile on ${\mathbb X}$, and the hypersurface $W$
is sufficiently general in that linear system. Set
$$
\eta=\overline{\eta}|_W\colon W\to S.
$$
The variety $W$ by the claim (i) of Theorem 0.1 has ay most
quadratic singularities of rank $\geqslant 8$, and for that reason
is a factorial variety with terminal singularities. Therefore,
$\eta\colon W\to S$ is a Fano-Mori fibre space, the fibres of
which are multiple projective spaces of index 1. Let $\eta'\colon
W'\to S'$ be an arbitrary rationally connected fibre space, that
is, a morphism of projective algebraic varieties, where the base
$S'$ and the fibre of general position $(\eta')^{-1}(s')$, $s'\in
S'$, are rationally connected, and moreover, $\mathop{\rm dim\,}
W'=\mathop{\rm dim\,} W$. Now \cite[Theorem 1.1]{Pukh15a},
combined with Theorem 0.1, immediately gives the following
result.\vspace{0.1cm}

{\bf Theorem 0.2.} {\it Assume that the Fano-Mori fibre space
$\eta\colon W\to S$ satisfies the following condition: for every
mobile family $\overline{\cal C}$ of curves on the base $S$,
sweeping out $S$, and a general curve $\overline{C}\in
\overline{\cal C}$ the class of an algebraic cycle
$$
-N (K_W\cdot \eta^{-1}(\overline{C}))-W_s
$$
is not effective, that is, it is not rationally equivalent to an
effective cycle of dimension $M$. Then every birational map
$\chi\colon W\dashrightarrow W'$ onto the total space of the
rationally connected fibre space $W'/S'$ (if such maps exist) is
fibre-wise, that is, there is a rational dominant map $\zeta\colon
S\dashrightarrow S'$, such that the following diagram commutes:}
$$
\begin{array}{rcccl}
   & W & \stackrel{\chi}{\dashrightarrow} & W' & \\
\eta\!\! & \downarrow &   &   \downarrow & \!\!\eta' \\
   & S & \stackrel{\zeta}{\dashrightarrow} & S'.
\end{array}
$$

{\bf Corollary 0.1.} {\it In the assumptions of Theorem 0.2 on the
variety $W$ there are no structures of a rationally connected
fibre space (and, the more so, of a Fano-Mori fibre space), the
fibre of which is of dimension less than $M$. In particular, the
variety $W$ is non-rational and every birational self-map of the
variety $W$ commutes with the projection $\eta$ and for that
reason induces a birational self-map of the base
$S$.}\vspace{0.1cm}

The condition for the cycles of dimension $M$, described in
Theorem 0.2, is satisfied if the linear system $|W|$ is
sufficiently mobile on ${\mathbb X}$. Let us demonstrate it by an
especially visual example, when ${\mathbb X}=\overline{{\mathbb
P}}\times S$ is the trivial fibre space over $S$. Let
$o^*=(0:\dots:0:1)=(0^{M+1}:1)\in\overline{\mathbb P}$ be the only
singular point of the weighted projective space $\overline{\mathbb
P}$. Consider the projection ``from the point $o^*$''
$$
\pi_{\mathbb P}\,\colon\overline{\mathbb P}\backslash\{o^*\}
\to{\mathbb P},
$$
where $\pi_{\mathbb P}((x_0:\dots:x_M:\xi))=(x_0:\dots:x_M)$. Let
$\overline{H}$ be the $\pi_{\mathbb P}$-pull back of the class of
a hyperplane in ${\mathbb P}$ on $\overline{{\mathbb P}}$. The
pull back of the class $\overline{H}$ on ${\mathbb
X}=\overline{{\mathbb P}}\times S$ with respect to the projection
onto the first factor we denote for simplicity by the same symbol
$\overline{H}$. Now
$$
\mathop{\rm Pic} {\mathbb X} = {\mathbb Z} \overline{H}\oplus
\overline{\eta}^* \mathop{\rm Pic} S,
$$
so that for some class $R\in \mathop{\rm Pic} S$ the relation
$$
W\sim dl\overline{H}+\overline{\eta}^* R
$$
holds and for that reason
$$
K_W=-\overline{H}|_W+\eta^* (R+K_S).
$$
This implies that the condition of Theorem 0.2 holds if for any
mobile family of curves $\overline{\cal C}$, sweeping out $S$, and
a general curve $\overline{C}\in \overline{\cal C}$, the
inequality
$$
((R+K_S)\cdot \overline{C})\geqslant 0
$$
holds. Therefore, the following claim is true.\vspace{0.1cm}

{\bf Theorem 0.3.} {\it Assume that the class $R+K_S\in
\mathop{\rm Pic} S$ is pseudo-effective and for every point $s\in
S$ we have $\eta^{-1}(s)=W_s\in {\cal F}_{\rm reg}$. Then in the
notations of Theorem 0.2 every birational map $\chi\colon W
\dashrightarrow W'$ is fibre-wise. In particular, every birational
self-map $\chi\in \mathop{\rm Bir} W$ induces a birational
self-map of the base $S$.}\vspace{0.1cm}

Another standard application of Theorem 0.1 is given by the
theorem on birational geometry of Fano direct products
\cite[Theorem 1]{Pukh05}. Recall that the following statement is
true.\vspace{0.1cm}

{\bf Theorem 0.4.} {\it Assume that primitive Fano varieties
$V_1,\dots, V_N$ satisfy the following properties:\vspace{0.1cm}

(i) for every effective divisor
$$
D_i\sim -nK_{V_i}
$$
the pair $(V_i,\frac{1}{n}D_i)$ is log canonical,\vspace{0.1cm}

(ii) for every mobile linear system
$$
\Sigma_i\subset |-nK_{V_i}|
$$
and a general divisor $D_i\in \Sigma_i$ the pair
$(V_i,\frac{1}{n}D_i)$ is canonical.\vspace{0.1cm}

Then on the direct product
$$
V_1\times\cdots\times V_N
$$
there are no other structures of a rationally connected fibre
space, apart from projecti\-ons on to direct fibres
$V_{i_1}\times\cdots\times V_{i_k}$.}\vspace{0.1cm}

The property (ii) was shown in \cite{Pukh19a} for a wider class of
multiple projective spaces than the one that is considered in this
paper. Of course, Theorem 0.1 implies that the conditions (i) and
(ii) are satisfied for every variety $V\in {\cal F}_{\rm reg}$.
Therefore, every variety considered in the present paper can be
taken as a factor of the direct product in Theorem 0.4.
\vspace{0.3cm}

%%%%%%%%%%%%%%%%%%%%%%%%%%%%%%%%%%%%%%%%%%%%%%%%%%%%%%%%
%%%%%%%%%%%%%%%%%%%%%%   subsection 0.2

{\bf 0.2. The regularity conditions.} The open subset ${\cal
F}_{\rm reg}$ are given by explicit local regularity conditions,
which we will now describe. To begin with, let us introduce an
auxiliary integral-valued parameter $\rho\in\{1,2,3,4\}$,
depending on $(d,l)$. Its meaning, the number reductions to a
hyperplane section, used in the proof of Theorem 0.1, will become
clear later. Set $\rho=4$, if $d=4$ and $21\leqslant l\leqslant
25$ and $\rho=1$, if $d\geqslant 18$ and $l\geqslant 2$. For the
remaining possible pairs $(d,l)$ the value $\rho\geqslant 2$ is
given by the following table:\vspace{0.5cm}

\begin{center}
\begin{tabular}{|c|c|c|}
\hline
$d$  & $l$ & $\rho$ \\
\hline 4 & $\geqslant 26$ & 3 \\
\hline 5 & 5, 6,\dots, 15 & 3 \\
\hline 5 & $\geqslant 16$ & 3\\
\hline 6 & 6 & 3\\
\hline 6 & $\geqslant 7$ & 2 \\
\hline
\end{tabular}
\end{center}
\vspace{0.5cm}

\begin{center}
\begin{tabular}{|c|c|c|}
\hline
$d$  & $l$ & $\rho$ \\
\hline 7 & 4 & 3 \\
\hline 7 & $\geqslant 5$ & 2 \\
\hline 8 & $\geqslant 4$ & 2 \\
\hline 9 & $\geqslant 3$ & 2 \\
\hline 10 & 3, 4, \dots, 17 & 2\\
\hline 11 & 2, 3, \dots, 8 & 2\\
\hline 12 & 2, 3, 4, 5 & 2\\
\hline 13 & 2, 3, 4 & 2\\
\hline 14 & 2, 3 & 2\\
\hline 15, 16, 17 & 2 & 2 \\
\hline
\end{tabular}
\end{center}
\vspace{0.5cm}

If the pair $(d,l)$ is not in the table, then $\rho=1$ (for
instance, for $d=14$, $l\geqslant 4$).\vspace{0.1cm}

One more table gives the function $\varepsilon(d,l)$, bounding
from below the codimension of the complement to the set ${\cal
F}_{\rm reg}$. Write this function as a function of the dimension
$M=(d-1)l$, for each of the possible values of the parameter
$\rho$ defined above.\vspace{0.3cm}

\begin{center}
\begin{tabular}{|c|c|}
\hline
$\rho$  & $\varepsilon(d,l)$ \\
\hline 1 & $\frac12(M^2-17M+56)$ \\
\hline 2 & $\frac12(M^2-21M+76)$ \\
\hline 3 & $\frac12(M^2-25M+90)$ \\
\hline 4 & $\frac12(M^2-31M+132)$ \\
\hline
\end{tabular}
\end{center}
\vspace{0.5cm}

Now let us state the regularity conditions.\vspace{0.1cm}

Let $o\in V$ be some point. The coordinate system $(x_0:x_1:\cdots
:x_M:\xi)$ can be chosen in such a way that
$$
o=(1:0:\cdots :0:0)
$$
(see \cite[\S 1]{Pukh19a}). The corresponding affine coordinates
are
$$
z_i=x_i/x_0,\quad i=1,\dots, M,\quad y=\xi/x^l.
$$
Now in the affine chart $\{x_0\neq 0\}={\mathbb
A}^{M+1}_{z_1,\dots,z_M,y}$ the hypersurface $V\cap\{x_0\neq 0\}$
is given by the equation $f=0$, where
$$
f=y^d+a_1(z_*)y^{d-1}+\dots+a_{d-1}(z_*)y+a_d(z_*),
$$
where the (non-homogeneous) polynomial $a_i(z_*)$ is of degree
$\leqslant il$. Furthermore, the following fact is true
(\cite[]{Pukh19a}): for any homogeneous polynomial
$\gamma(x_0,\dots,x_M)$ of degree $l$ the equation
$\xi=\gamma(x_*)$ defines a hypersurface
$R_{\gamma}\subset\overline{\mathbb P}$ that does not contain the
point $o^*=(0^{M+1}:1)$, and moreover the projection
$$
\pi_{\mathbb P}|_{R_{\gamma}}\colon R_{\gamma}\to {\mathbb P}
$$
is an isomorphism. In this way, the hypersurface $V\cap
R_{\gamma}$ in $R_{\gamma}$ identifies naturally with a
hypersurface in ${\mathbb P}={\mathbb P}^M$, and its intersection
with the affine chart $\{x_0\neq 0\}$ identifies with a
hypersurface in the affine space ${\mathbb A}^M_{z_1,\dots,z_M}$.
The regularity conditions, given below, are assumed to be
satisfied for the hypersurface $V_{\gamma}=V\cap R_{\gamma}$ for a
general polynomial $\gamma(x_*)$.\vspace{0.1cm}

Assume that the point $o\in V$ is non-singular, so that $o\in
V_{\gamma}$ is non-singular, too. Let $P\subset {\mathbb A}^M$ be
an arbitrary linear subspace of codimension $(\rho - 1)$, that is,
$o\in P$, that is not contained in the tangent hyperplane
$T_oV_{\gamma}$. Let
$$
f_P=q_1+q_2+\cdots +q_{dl}
$$
be the affine equation of the hypersurface $P\cap V_{\gamma}$,
which is non-singular at the point $o$, decomposed into
homogeneous components (with respect to an arbitrary system of
linear coordinates on $P$).\vspace{0.1cm}

(R1.1) For any linear form
$$
\lambda\not\in\langle q_1\rangle
$$
the sequence of homogeneous polynomials
$$
q_1|_{\{\lambda=0\}},\quad q_2|_{\{\lambda=0\}},\quad \dots,\quad
q_{M-\rho-2}|_{\{\lambda=0\}}
$$
is regular in the local ring ${\cal O}_{o,P}$.\vspace{0.1cm}

(R1.2) The linear span of any irreducible component of the closed
set
$$
\{q_1=q_2=q_3=0\}
$$
is the hyperplane $\{q_1=0\}$.\vspace{0.1cm}

(R1.3) For any linear form $\lambda\not\in \langle q_1\rangle$ the
set
$$
\overline{P\cap V_{\gamma}\cap\{q_1=q_2=0\}\cap\{\lambda=0\}}
$$
is irreducible and reduced.\vspace{0.1cm}

(R1.4) If $\rho\geqslant 2$, then the rank of the quadratic form
$$
q_2|_{\{q_1=0\}}
$$
is at least $8+2(\rho-2)$.\vspace{0.1cm}

We say that a non-singular point $o\in V$ is {\it regular}, if for
a general polynomial $\gamma(x_*)$ and any subspace $P\not\subset
T_oV_{\gamma}$ the conditions (R1.1-3) are
satisfied.\vspace{0.1cm}

Assume now that the point $o\in V$ is singular, so that the
hypersurface $V_{\gamma}$ is also singular at that
point.\vspace{0.1cm}

(R2.1) The point $o\in V_{\gamma}$ is a quadratic singularity of
rank $\geqslant 2\rho+6$.\vspace{0.1cm}

Let $P\subset {\mathbb A}^M$ be an arbitrary linear subspace of
codimension $\rho+2$, that is, $o\in P$, and
$$
f_P=q_2+q_3+\cdots +q_{dl}
$$
is the affine equation of the hypersurface $P\cap V_{\gamma}$,
decomposed into homogeneous components (in particular, $q_2$ is a
quadratic form of rank $\geqslant 2$).\vspace{0.1cm}

(R2.2) The sequence of homogeneous polynomials
$$
q_2,\quad q_3,\quad\dots,\quad q_{M-\rho-4}
$$
is regular in the local ring ${\cal O}_{o,P}$.\vspace{0.1cm}

We say that a singular point $o\in V$ is {\it regular}, if for a
general polynomial $\gamma(x_*)$ and any subspace $P\subset
{\mathbb A}^M$ of codimension $\rho+2$ the conditions (R2.1,2)
hold.\vspace{0.1cm}

Finally, we say that the variety $V$ is {\it regular}, if it is
regular at every point $o\in V$, singular or non-singular. Set
$$
{\cal F}_{\rm reg}\subset {\cal F}
$$
to be the Zariski open subset of regular hypersurfaces (that it is
non-empty, follows from the estimate for the codimension of the
complement). Obviously, every hypersur\-face $V\in {\cal F}_{\rm
reg}$ has at worst quadratic singularities of rank $\geqslant 8$,
so that the claim (i) of Theorem 0.1 is true.\vspace{0.3cm}

%%%%%%%%%%%%%%%%%%%%%%%%%%%%%%%%%%%%%%%%%%%%%%%%%%%%%%%%%%%%%%%%%
%%%%%%%%%%%%%%%%%%%%%   subsection 0.3

{\bf 0.3. The structure of the paper, historical remarks and
acknowledge\-ments.} A proof of the claim (ii) of Theorem 0.1 is
given in Subsections 1.2 and 1.3. A proof of the claim (iii) of
Theorem 0.1 in Subsection 1.1 is reduced to two facts about
hypersurfaces in the projective space ${\mathbb P}^N$, which are
applied to the hypersurface $V_{\gamma}\subset {\mathbb P}$, both
in the singular and non-singular cases. Proofs of those two facts
are given, respectively, in \S 2 and \S 3.\vspace{0.1cm}

The equality of the global (log) canonical threshold to one is
shown for many families of primitive Fano varieties, starting from
the pioneer paper \cite{Pukh05} (for a general variety in the
family). For Fano complete intersections in the projective space
the best progress in that direction (in the sense of covering the
largest class of families) was made in \cite{Pukh18a}. The double
covers were considered in \cite{Pukh08a}. Fano three-folds,
singular and non-singular, were studied in the papers
\cite{Ch2009b,ChShr2008,ChParkWon2014,Ch08} and many others.
However, the non-cyclic covers of index 1 in the arbitrary
dimension were never studied up to now: the reason, as it was
explained in \cite{Pukh19a}, was that the technique of
hypertangent divisors does not apply to these varieties in a
straightforward way. As it turned out (see \cite{Pukh19a}), the
technique of hypertangent divisors should be applied to a certain
subvariety, which identifies naturally with a hypersurface (of
general type) in the projective space. This approach is used in
the present paper, too.\vspace{0.1cm}

The author thanks The Leverhulme Trust for the support of the
present work (Research Project Grant RPG-2016-279).\vspace{0.1cm}

The author is also grateful to the colleagues in the Divisions of
Algebraic Geometry and Algebra at Steklov Institute of Mathematics
for the interest to his work, and to the colleagues-algebraic
geometers at the University of Liverpool for the general support.

%%%%%%%%%%%%%%%%%%%%%%%%%%%%%%%%%%%%%%%%%%%%%%%%%%%%%%%%%%%%%%%%%
%%%%%%%%%%%%%%%%%%%%%%%%%%%%%%%%%%%%%%%%%%%%%%%%%%%%%%%%%%%%%%%%%
%%%%%%%%%%%%%%%%%%%%   SECTION 1

\section{Proof of the main result}

In Subsection 1.1 the proof of part (iii) of Theorem 0.1 is
reduced to two intermediate claims, the proofs of which are given
in \S 2 and \S 3. In Subsections 1.2,3 we show part (ii) of
Theorem 0.1. First (Subsection 1.2) we give the estimates for the
codimension of the sets of polynomials, violating each of the
regularity conditions, after that (Subsection 1.3) we explain how
to obtain these estimates.\vspace{0.3cm}

{\bf 1.1. Exclusion of maximal singularities.} Fix the parameters
$d,l$. Recall that the integer $\rho\in\{1,2,3,4\}$ depends on
$d,l$ (see the table in Subsection 0.2). Fix a variety $V\in{\cal
F}_{\rm reg}$. Assume that $D\sim nH$ is an effective divisor on
$V$, such that the pair $(V,\frac{1}{n}D)$ is not canonical. Our
aim is to get a contradiction. This would prove the claim
(iii).\vspace{0.1cm}

If there is a non-canonical singularity of the pair
$(V,\frac{1}{n}D)$, the centre of which is of positive dimension,
then the pair
$$
(\Gamma,\frac{1}{n}D_{\Gamma}),
$$
where $\Gamma=V_{\gamma}$ for some polynomial $\gamma(x_*)$ of
general position (see Subsection 0.2) and
$D_{\Gamma}=D|_{\Gamma}$, is again non-canonical. If the centres
of all non-canonical singularities of the pair $(V,\frac{1}{n}D)$
are points, let us take a polynomial $\gamma(x_*)$ of general
position such that the hypersurface $\Gamma=V_{\gamma}$ contains
one of them. In that case the pair
$(\Gamma,\frac{1}{n}D_{\Gamma})$ is even non log
canonical.\vspace{0.1cm}

In any case we obtain a factorial hypersurface
$\Gamma\subset{\mathbb P} ={\mathbb P}^M$ of degree $dl$ with at
worst quadratic singularities of rank $\geqslant 2\rho+6\geqslant
8$, and an effective divisor $D_{\Gamma}\sim nH_{\Gamma}$ on it
(where $H_{\Gamma}$ is the class of a hyperplane section, so that
$\mathop{\rm Pic}\Gamma={\mathbb Z}H_{\Gamma}$), such that the
pair $(\Gamma,\frac{1}{n}D_{\Gamma})$ is non-canonical. Now we
work only with that pair, forgetting about the original variety
$V$ (within the limits of the proof of the claim (iii) of Theorem
0.1). Let
$$
\mathop{\rm CS}\left(\Gamma,\frac{1}{n}D_{\Gamma}\right)
$$
be the union of the centres of all non-canonical singularities of
that pair.\vspace{0.1cm}

{\bf Proposition 1.1.} {\it The closed set $\mathop{\rm
CS}(\Gamma,\frac{1}{n}D)$ is contained in the singular locus
$\mathop{\rm Sing}\Gamma$ of the hypersurface
$\Gamma$.}\vspace{0.1cm}

{\bf Proof} makes the contents of \S 2.\vspace{0.1cm}

Therefore,
$$
\mathop{\rm codim}\left(\mathop{\rm
CS}\left(\Gamma,\frac{1}{n}D_{\Gamma}\right)\subset
\Gamma\right)\geqslant 7.
$$
Let us define a sequence of rational numbers $\alpha_k$,
$k\in{\mathbb Z}_+$, in the following way:
$$
\alpha_0=1,\quad \alpha_{k+1}=\frac12\alpha_k+1.
$$
(We can simply write $\alpha_k=2-\frac{1}{2^k}$, but for us it is
important how $\alpha_{k+1}$ and $\alpha_k$ are related.) In order
to exclude the maximal (non-canonical) singularities, we will need
only the four values:
$$
\alpha_1=\frac32,\,\,\alpha_2=\frac74,\,\,\alpha_3=\frac{15}{8},
\,\,\alpha_4=\frac{31}{16}.
$$
Let $o\in\Gamma$ be a point of general position on the irreducible
component of maximal dimension of the closed set $\mathop{\rm
CS}(\Gamma,\frac{1}{n}D_{\Gamma})$. Consider a general
5-dimensional subspace in ${\mathbb P}$, containing the point $o$.
Let $P$ be the section of the hypersurface $\Gamma$ by that
subspace. Obviously, $P\subset{\mathbb P}^5$ is a hypersurface of
degree $dl$ with a unique singular point, a non-degenerate
quadratic point $o$. Denoting $D_{\Gamma|_P}$ by the symbol $D_P$,
we get $D_P\sim nH_P$, where $H_P$ is the class of a hyperplane
section. By the inversion of adjunction, the point $o$ is the
centre of a non {\it log} canonical singularity of the pair
$(P,\frac{1}{n}D_P)$, and moreover,
$$
\mathop{\rm LCS}\left(P,\frac{1}{n}D_P\right)=\{o\}.
$$
This implies that
$$
\mathop{\rm mult}\nolimits_o D_P>2n
$$
and therefore
$$
\mathop{\rm mult}\nolimits_o D_{\Gamma}>2n=2\alpha_0n.
$$

{\bf Proposition 1.2.} {\it There is a sequence of irreducible
varieties $\Gamma_i$, $i=0,1,\dots,\rho$, such that:\vspace{0.1cm}

{\rm (i)} $\Gamma_0=\Gamma$ and $\Gamma_{i+1}$ is a hyperplane
section of the hypersurface $\Gamma_i\subset{\mathbb P}^{M-i}$,
containing the point $o$,\vspace{0.1cm}

{\rm (ii)} on the variety $\Gamma_{\rho}$ there is a prime divisor
$D^*\sim n^*H^*$, where $H^*$ is the class of a hyperplane section
of the hypersurface $\Gamma_{\rho}$, satisfying the inequality}
$$
\mathop{\rm mult}\nolimits_o D^*>2\alpha_{\rho}n^*.
$$

{\bf Proof} makes the contents of \S 3.\vspace{0.1cm}

Note that by the condition (R2.1) all hypersurfaces
$\Gamma_1,\dots,\Gamma_{\rho}$ are factorial, so that $\mathop{\rm
Pic}\Gamma_{\rho}={\mathbb Z}H^*$. Furthermore, $\rho\geqslant 1$,
so that
$$
\mathop{\rm mult}\nolimits_o D^*>3n.
$$
Now let us consider general hypertangent divisors
$D_2,\dots,D_{M-\rho-2}$ on the hypersurface $\Gamma_{\rho}$ (for
the definition and construction of hypertangent divisors, see
\cite[Chapter 3]{Pukh13a}) and construct in the usual way a
sequence of irreducible subvarieties $Y_i\subset \Gamma_{\rho}$ of
codimension $i=1,2,\dots,M-\rho-3$, such that $Y_1=D^*$, the
subvariety $Y_2$ is an irreducible component of the effective
cycle $(Y_1\circ D_2)$ with the maximal ratio of the multiplicity
$\mathop{\rm mult}_o$ to the degree $\mathop{\rm deg}$, and for
$i=3,\dots,M-\rho-3$ the subvariety $Y_i$ is an irreducible
component of the effective cycle $(Y_{i-1}\circ D_{i+1})$ with the
maximal value of the ratio of the multiplicity $\mathop{\rm
mult}_o$ to the degree. That it is possible to go through with
this construction, is ensured by the condition (R2.2). Note that
the first step of this construction is possible because the
hypertangent divisor $D_2$ is irreducible, $D_2\sim 2H^*$ and the
equality
$$
\mathop{\rm mult}\nolimits_oD_2=6=3\cdot 2
$$
holds, so that $Y_1\neq D_2$. The hypertangent divisor $D_3$ does
not take part in the construction.\vspace{0.1cm}

For the irreducible surface
$$
S=Y_{M-\rho-3}\subset\Gamma_{\rho}
$$
we get the estimate
$$
\frac{\mathop{\rm mult}_o}{\mathop{\rm deg}}S>\frac{1}{dl}\cdot
2\alpha_{\rho}\cdot\frac32\cdot\frac54\cdot\frac65\cdot\dots\cdot
\frac{M-\rho-1}{M-\rho-2}=
\frac{3(M-\rho-1)}{4dl}\,\alpha_{\rho}\geqslant 1,
$$
which is impossible (the last inequality checks directly for each
of the possible values of $\rho$ and the corresponding values of
$d,l$). Thus we obtained a contradiction, which completes the
proof of the claim (iii) of Theorem 0.1.\vspace{0.3cm}

%%%%%%%%%%%%%%%%%%%%%%%%%%%%%%%%%%%%%%%%%%%%%%%%%%%%%%%%%%%%%%
%%%%%%%%%%%%%%%%%%%%%   subsection 1.2

{\bf 1.2. Estimating the codimension of the set ${\cal
F}\backslash{\cal F}_{\rm reg}$.} Let us prove the claim (ii) of
Theorem 0.1. Denote by the symbol ${\cal F}_{i.j}$ the closure of
the set of hypersurfaces $V\in{\cal F}$, violating the condition
(Ri.j) at at least one point. Here
$$
i.j\in\{1.1,1.2,1.3,1.4,2.1,2.2\}.
$$
For these values of $i.j$ we set, respectively,
$$
\varepsilon_{i.j}=\mathop{\rm codim}({\cal F}_{i.j}\subset{\cal
F}).
$$
We omit the symbols $d,l$ in order to simplify the formulas,
however $\varepsilon_{i.j}=\varepsilon_{i.j}(d,l)$ are functions
of these parameters. The following claim is true.\vspace{0.1cm}

{\bf Proposition 1.3.} {\it The following inequalities
hold:}\vspace{0.1cm}

(i)\,\,$\varepsilon_{1.1}\geqslant
\frac12(M^2-(4\rho+5)M+(3\rho^2+3\rho))$,\vspace{0.1cm}

(ii)\,\,$\varepsilon_{1.2}\geqslant
\frac12(M^2-(4\rho+11)M+(3\rho^2-15\rho+32))$,\vspace{0.1cm}

(iii)\,\,$\varepsilon_{1.3}\geqslant
\frac12(M^2-(4\rho+13)M+(3\rho^2+11\rho+42))$,\vspace{0.1cm}

(iv)\,\,$\varepsilon_{1.4}\geqslant
\frac12(M^2-(4\rho+9)M+(4\rho^2+14\rho+16))$,\vspace{0.1cm}

(v)\,\,$\varepsilon_{2.1}\geqslant
\frac12(M^2-(6\rho+7)M+(4\rho^2+14\rho+12))$,\vspace{0.1cm}

(vi)\,\,$\varepsilon_{2.2}\geqslant
\frac12(M^2-(4\rho+1)M+(3\rho^2-\rho))$.\vspace{0.1cm}

{\bf Proof.} The regularity conditions must be satisfied for {\it
any} point $o$, {\it any} linear subspace $P$ of the required
codimension and {\it any} linear form $\lambda$ (the polynomial
$\gamma(x_*$) is assumed to be general and does not influence the
estimating of the codimension of the sets ${\cal F}_{i.j}$).
Therefore, the problem of getting a lower bound for the numbers
$\varepsilon_{i.j}$ reduces obviously to a similar problem for
varieties $V\in{\cal F}$ violating the condition (Ri.j) at a {\it
fixed} point $o$, for a {\it fixed} linear subspace and a {\it
fixed} linear form $\lambda$. The solution of the latter problem
comes from the claims of Propositions 1.4 and 1.5, shown below.
More precisely, the estimates for the conditions (R1.4) and (R2.1)
follow from the part (i) of Proposition 1.4, for the condition
(R1.2) from the part (ii) of Proposition 1.4, for the condition
(R1.3) from the part (iii) of that proposition. The estimates for
the conditions (R1.1) and (R2.2) follow from Proposition 1.5. The
proof is complete. Q.E.D.\vspace{0.1cm}

Now, in order to prove the claim (ii) of Theorem 0.1, it is
sufficient to check that the function $\varepsilon(d,l)$ is the
minimum of the right hand sides in the inequalities (i)-(vi) of
Proposition 1.3. This work is elementary and we do not give it
here. Q.E.D. for the claim (ii) of Theorem 1.\vspace{0.3cm}

%%%%%%%%%%%%%%%%%%%%%%%%%%%%%%%%%%%%%%%%%%%%%%%%%%%%%%%%%%%%%%%%%%%
%%%%%%%%%%%%%%%%%%%%%%%%   subsection 1.3

{\bf 1.3. Quadratic forms and regular sequences.} By the symbol
${\cal P}_{i,N}$ we denote the linear space of homogeneous
polynomials of degree $i\in{\mathbb Z}_+$ in $N$ variables
$u_1,\dots,u_N$. For $i\leqslant j$ we write
$$
{\cal P}_{[i,j],N}=\bigoplus^j_{k=i}{\cal P}_{k,N},
$$
and ${\cal P}_{\leqslant i,N}=\bigoplus\limits^i_{k=0}{\cal
P}_{k,N}$. The number of variables $N$ is fixed, so we omit the
symbol $N$ and write ${\cal P}_k,{\cal P}_{[i,j]}$ and so on. Let
$$
{\cal X}_{2,\leqslant r}\subset{\cal P}_2
$$
be the closed subset of quadratic forms of rank $\leqslant r$. Let
$$
{\cal X}_{2,3}\subset{\cal P}_{[2,3]}
$$
be the closed subset of pairs $(w_2,w_3)$, such that the closed
set $\{w_2=w_3=0\}\subset{\mathbb P}^{N-1}$ has at least one
degenerate component (that is, a component, the linear span of
which is of dimension $\leqslant N-2$). Let $Q\subset{\mathbb
P}^{N-1}$ be a factorial quadric. For $m\geqslant 4$ let
$$
{\cal X}_{m,Q}\subset{\cal P}_m
$$
be the closed subset of polynomials $w_m$, such that the divisor
$\{w_m|_Q=0\}$ on $Q$ is reducible or non-reduced. The following
claim is true.\vspace{0.1cm}

{\bf Proposition 1.4.} (i) {\it The following equality holds:}
$$
\mathop{\rm codim}({\cal X}_{2,\leqslant r}\subset{\cal P}_2)=
{N-r+1\choose 2}.
$$

(ii) {\it The following inequality holds:}
$$
\mathop{\rm codim}({\cal X}_{2,3}\subset{\cal P}_{[2,3]})\geqslant
{N-3\choose 2}.
$$

(iii) {\it The following inequality holds:}
$$
\mathop{\rm codim}({\cal X}_{m,Q}\subset{\cal P}_m)\geqslant 2
{m+N-4\choose N-3}.
$$

{\bf Proof.} The claim (i) is well known. Let us show the
inequality (ii). Taking into account the part (i), we may assume
that the quadratic form $w_2$ is of rank $\geqslant 5$, so that
the quadric $\{w_2=0\}$ is factorial. If the closed set
$w_2=w_3=0$ has a degenerate component, then the divisor
$\{w_3|_{\{w_2=0\}}=0\}$ on the quadric $\{w_2=0\}$ is either
reducible, or non-reduced, so that in any case it is a sum of a
hyperplane section and a section of the quadric $\{w_2=0\}$ by
some quadratic hypersurface. Calculating the dimensions of the
corresponding linear systems, we get that for a fixed quadratic
form $w_2$ of rank $\geqslant 5$ the closed set of polynomials
$w_3\in{\cal P}_3$, such that the divisor $\{w_3|_{\{w_2=0\}}=0\}$
is reducible or non-reduced, is of codimension
$$
{N+2\choose 3}-{N+1\choose 2}-2N+2
$$
in ${\cal P}_3$. It is easy to see that this expression is higher
than the right hand side of the inequality (ii). This proves the
claim (ii).\vspace{0.1cm}

Let us show the inequality (iii). Recall that the quadric
$Q\subset{\cal P}^{N-1}$ is assumed to be factorial (that is, the
rank of the corresponding quadratic form is at least 5). Set
$h_Q(m)=h^0(Q,{\cal O}_Q(m))$ for $m\geqslant 1$. It is easy to
check that
$$
h_Q(m)=\frac{(m+(N-3))\dots(m+1)}{(N-2)!}(2m+(N-2))
$$
is a polynomial in $m$ with positive coefficients. This implies
that for $0<s<t\leqslant\frac12m$ the inequality
$$
h_Q(t)-h_Q(s)<h_Q(m-s)-h_Q(m-t)
$$
holds, which can be re-written as
$$
h_Q(t)+h_Q(m-t)<h_Q(s)+h_Q(m-s).
$$
If the divisor $\{w_m|_Q=0\}$ is not irreducible and reduced, then
it is a sum of two effective divisors on $Q$, which are cut out on
$Q$ by hypersurfaces of degree $1\leqslant a\leqslant\frac12m$ and
$(m-a$). For that reason,
$$
\mathop{\rm dim}{\cal X}_{m,Q}=\mathop{\rm max}_{1\leqslant
a\leqslant\frac12m}\{h_Q(a)+h_Q(m-a)\}.
$$
By what was said above, the right hand side of that inequality is
$h_Q(1)+h_Q(m-1)$, so that
$$
\mathop{\rm codim}({\cal X}_{m,Q}\subset{\cal
P}_m)=h_Q(m)-h_Q(m-1)-h_Q(1).
$$
Elementary computations show that the right hand side of the last
equality is
$$
\frac{(m+(N-4))\dots(m+1)}{(N-3)!}(2m+(N-3))-N,
$$
which is certainly higher than
$$
2{m+N-4\choose N-3}.
$$
Proof of Proposition 1.4 is complete. Q.E.D.\vspace{0.1cm}

{\bf Remark 1.1.} The estimates for
$\varepsilon_{1.2},\varepsilon_{1.4}$ and $\varepsilon_{2.1}$,
given in Proposition 1.3, are obtained from the claims (i) and
(ii) of Proposition 1.4 by elementary computations. It is slightly
less obvious, how to obtain the estimate for $\varepsilon_{1.3}$,
starting from the claim (iii) of Proposition 1.4, for that reason
we will explain briefly, how to do it. Fixing the linear subspace
$P$ and the linear forms $q_1$ and $\lambda$, consider the quadric
\begin{equation}\label{29.05.2019.1}
q_2|_{P\cap\{q_1=\lambda=0\}}=0.
\end{equation}
The codimension of the set of quadratic forms, for which this
quadric is of rank $\leqslant 4$ and so not factorial, is given by
the claim (i) of Proposition 1.4. It is from here that we get the
estimate for $\varepsilon_{1.3}$ in Proposition 1.3. It remains to
show that the violation of the condition (R1.3) {\it under the
assumption that the quadric (\ref{29.05.2019.1}) is factorial},
gives at least the same (in fact, much higher) codimension. It is
to the factorial quadric (\ref{29.05.2019.1}) that we apply the
estimate (iii) of Proposition 1.4. There is, however, a delicate
point here. The hypersurface $P\cap V_{\gamma}$ is given by a
polynomial that has at the point $o$ the linear part $q_1$ and the
quadratic part $q_2$, which both vanish when restricted onto the
quadric (\ref{29.05.2019.1}). The other homogeneous components
$q_3,\dots,q_{dl}$ are arbitrary. In the inequality (iii) of
Proposition 1.4 the codimension of the ``bad'' set ${\cal
X}_{m,Q}$ is considered with respect to the whole space ${\cal
P}_m$, whereas in order to prove the inequality (iii) of
Proposition 1.3, we need the codimension with respect to the space
of homogeneous polynomials of degree $dl$, the non-homogeneous
presentation of which at the fixed point $o$ has zero linear and
quadratic components. However, this does not make any influence on
the final result, because the codimension of the set ${\cal
X}_{m,Q}$ in ${\cal P}_m$ is very high.\vspace{0.1cm}

Now let for $2\leqslant k\leqslant N-2$
$$
{\cal X}_{[2,k]}\subset{\cal P}_{[2,k]}
$$
be the set of non-regular tuples $(h_2,\dots,h_k)$ of length
$k-1\leqslant N-3$, where $h_i\in{\cal P}_i={\cal P}_{i,N}$, that
is, the system of equations
$$
h_2=\dots=h_k=0
$$
defines in ${\cal P}^{N-1}$ a closed subset of codimension
$\leqslant k-2$.\vspace{0.1cm}

{\bf Proposition 1.5.} {\it The following equality holds:}
$$
\mathop{\rm codim}({\cal X}_{[2,k]}\subset {\cal
P}_{[2,k]})={N+1\choose 2}
$$

{\bf Proof:} see \cite[Chapter 3, Section 1]{Pukh13a}.

%%%%%%%%%%%%%%%%%%%%%%%%%%%%%%%%%%%%%%%%%%%%%%%%%%%%%%%%%%%%%%%%%
%%%%%%%%%%%%%%%%%%%%%%%%%%%%%%%%%%%%%%%%%%%%%%%%%%%%%%%%%%%%%%%%%
%%%%%%%%%%%%%%%%%%%%   SECTION 2

\section{Exclusion of maximal singularities \\ at smooth points}

In this section we consider factorial hypersurfaces
$X\subset{\mathbb P}^N$, satisfying certain additional conditions.
We show that the centre of every non-canonical singularity of the
pair $(X,\frac{1}{n} D_X)$, where $D_X\sim nH_X$ is cut out on $X$
by a hypersurface of degree $n\geqslant 1$, is contained in the
singular locus $\mathop{\rm Sing} X$. In Subsection 2.1 we list
the conditions that are satisfied by the hypersurface $X$, state
the main result and exclude non-canonical singularities with the
centre of a small ($\leqslant 3$) codimension on $X$. In
Subsections 2.2 and 2.3, following (with minor modification) the
arguments of Subsection 2.1 in \cite{Pukh05}, we exclude
non-canonical singularities of the pair $(X,\frac{1}{n} D_X)$, the
centre of which is not contained in $\mathop{\rm Sing} X$. In
subsection 2.3 we use, for this purpose, the standard technique of
hypertangent divisors. As a first application, we obtain a proof
of Proposition 1.1.\vspace{0.3cm}

{\bf 2.1. Regular hypersurfaces.} Let $X\subset{\mathbb P}^N$,
where $N\geqslant 8$, be a hypersurface, satisfying the condition
$$
\mathop{\rm codim}(\mathop{\rm Sing}X\subset X)\geqslant 5.
$$
In particular, $X$ is factorial and $\mathop{\rm Pic}X={\mathbb
Z}H_X$, where $H_X$ is the class of a hyperplane section. Let
$o\in X$ be a non-singular point and
$$
z_1,\dots,z_N
$$
a system of affine coordinates on ${\mathbb A}^N\subset{\mathbb
P}^N$ with the origin at the point $o$, and the hypersurface $X$
in this coordinate system is given by the equation $h=0$, where
$$
h=h_1+h_2+\dots+h_{{\rm deg}X}
$$
and the polynomials $h_i$ are homogeneous of degree $i$. We assume
that the inequality
\begin{equation}\label{19.04.2019.5}
N-2\leqslant\mathop{\rm deg}X\leqslant\frac32(N-3)
\end{equation}
holds.\vspace{0.1cm}

Now let us state the regularity conditions for the hypersurface
$X$ at the point $o$.\vspace{0.1cm}

(N1) For any linear form
$$
\lambda(z_*)\not\in\langle h_1\rangle
$$
the sequence of homogeneous polynomials
$$
h_1|_{\{\lambda=0\}},\, h_2|_{\{\lambda=0\}},\,\dots,\,
h_{N-3}|_{\{\lambda=0\}}
$$
is regular (in the local ring ${\cal O}_{o,{\mathbb
P}^N}$).\vspace{0.1cm}

(N2) The linear span of every irreducible component of the closed
set
$$
h_1=h_2=h_3=0
$$
is the hyperplane $\{h_1=0\}$.\vspace{0.1cm}

(N3) For any linear form $\lambda\not\in\langle h_1\rangle$ the
set
$$ \overline{X\cap\{h_1=h_2=0\}\cap\{\lambda=0\}}
$$
is irreducible and reduced.\vspace{0.1cm}

{\bf Proposition 2.1.} {\it Assume that the hypersurface $X$
satisfied the conditions (N1-3) at every non-singular point $o\in
X$. Then for every pair $(X,\frac{1}{n}D_X)$, where $D_X\sim nH_X$
is an effective divisor, the union of the centres of all
non-canonical singularities $\mathop{\rm CS\,}(X,\frac{1}{n}D_X)$
of that pair is contained in the closed set $\mathop{\rm Sing}
X$.}\vspace{0.1cm}

{\bf Proof.} Assume the converse: for some effective divisor
$D_X\sim nH_X$
$$
\mathop{\rm CS\,}\left(X,\frac{1}{n}D_X\right)
\not\subset\mathop{\rm Sing} X.
$$
Let $Y$ be an irreducible component of the set $\mathop{\rm CS\,}
(X,\frac{1}{n}D_X)$, which is not contained in $\mathop{\rm Sing}
X$, the dimension of which is maximal among all such
components.\vspace{0.1cm}

{\bf Lemma 2.1.} {\it The following inequality holds:}
$$
\mathop{\rm codim}(Y\subset X)\geqslant 4.
$$

{\bf Proof.} Assume the converse:  $\mathop{\rm codim}(Y\subset
X)\leqslant 3$. Since $Y$ is the centre of some non canonical
singularity of the pair $(X,\frac{1}{n}D_X)$ and
$Y\not\subset\mathop{\rm Sing}X$, we get the inequality
$\mathop{\rm mult}\nolimits_YD_X>n$. Since the codimension of the
set $\mathop{\rm Sing} X$ is at least 5, we can take a curve
$C\subset X$, such that
$$
C\subset X\backslash\mathop{\rm Sing} X.
$$
Obviously, $\mathop{\rm mult}\nolimits_C D_X>n$. Now repeating the
arguments in the proof of Lemma 2.1 in \cite[Chapter 2]{Pukh13a}
word for word, we get a contradiction which completes the proof of
Lemma 2.1.\vspace{0.3cm}

%%%%%%%%%%%%%%%%%%%%%%%%%%%%%%%%%%%%%%%%%%%%%%%%%%%%%%%%%%%%%%%%%%
%%%%%%%%%%%%%%%%%%%%%%%%% subsection 2

{\bf 2.2. Restriction onto a hyperplane section.} Let $o\in Y$ be
a point of general position, $o\not\in\mathop{\rm Sing}X$.
Consider the section $P\subset X$ by a general linear subspace of
dimension 4, containing the point $o$. The hypersurface
$P\subset{\mathbb P}^4$ is non-singular, so that $\mathop{\rm
Pic}P={\mathbb Z}H_P$ by the Lefschetz theorem, where $H_P$ is the
class of a hyperplane section of the variety $P$. Set
$D_P=D_X|_P$, so that $D_P\sim nH_P$. By inversion of adjunction,
the pair $(P,\frac{1}{n}D_P)$ is not log canonical; moreover, by
construction,
$$
\mathop{\rm LCS\,}\left(P,\frac{1}{n}D_P\right)=\{o\}.
$$
Let $\varphi_P\colon P^+\to P$ be the blow up of the point $o$,
$E_P=\varphi^{-1}_P(o)\cong{\mathbb P}^2$ the exceptional divisor,
$D^+_P$ the strict transform of the divisor $D_P$ on
$P^+$.\vspace{0.1cm}

{\bf Lemma 2.2.} {\it There is a line $L\subset E_P$, satisfying
the inequality}
$$
\mathop{\rm mult}\nolimits_oD_P+\mathop{\rm
mult}\nolimits_LD_P^+>2n.
$$

{\bf Proof.} This follows from \cite[Proposition 9]{Pukh05}.
Q.E.D.\vspace{0.1cm}

The blow up $\varphi_P$ can be viewed as the restriction onto the
subvariety $P$ of the blow up $\varphi_X\colon X^+\to X$ of the
point $o$ with the exceptional divisor $E_X\cong{\mathbb
P}^{N-2}$. Lemma 2.2 implies that there is a hyperplane
$\Theta\subset E_X$, satisfying the inequality
\begin{equation}\label{18.04.2019.3}
\mathop{\rm mult}\nolimits_oD_X+ \mathop{\rm
mult}\nolimits_{\Theta}D^+_X>2n.
\end{equation}
The rest of the proof of Proposition 2.1 repeats the proof of part
(i) of Theorem 2 in \cite[п. 2.1]{Pukh05} almost word for word.
For the convenience of the reader we briefly reproduce those
arguments. By the symbol $|H_X-\Theta|$ we denote the pencil of
hyperplane sections $R$ of the hypersurface $X$, such that $R\ni
o$ and $R^+\cap E_X=\Theta$ (where $R^+\subset X^+$ is the strict
transform). Let $R\in|H_X-\Theta|$ be a general element of the
pencil. Set $D_R=D_X|_R$.\vspace{0.1cm}

{\bf Lemma 2.3.} {\it The following inequality holds:}
\begin{equation}\label{19.04.2019.1}
\mathop{\rm mult}\nolimits_oD_R>2n.
\end{equation}

{\bf Proof.} This is Lemma 3 in \cite{Pukh05} (our claim follows
directly from the inequality (\ref{18.04.2019.3}) and the choice
of the section $R$). Q.E.D. for the lemma.\vspace{0.1cm}

Consider the tangent hyperplane $T_oR\subset{\mathbb P}^{N-1}$ to
the hypersurface $R$ at the point $o$. The intersection $T_R=R\cap
T_oR$ is a hyperplane section of $R$. Therefore, $T_R\sim H_R$ is
a prime divisor on $R$. By the condition (N1) the equality
$\mathop{\rm mult}_oT_R=2$ holds. Therefore, if
$$
D_R=aT_R+D^{\sharp}_R,
$$
where $a\in{\mathbb Z}_+$ and the effective divisor
$D^{\sharp}_R\sim(n-a)H_R$ does not contain $T_R$ as a component,
then the inequality
$$
\mathop{\rm mult}\nolimits_oD_R^{\sharp}>2(n-a)
$$
holds. In order not to make the notations too complicated, we
assume that $a=0$, that is, $D_R\sim nH_R$ does not contain $T_R$
as a component. Moreover, by the linearity of the inequality
(\ref{19.04.2019.1}) in $D_R$, we may assume that $D_R$ is a prime
divisor.\vspace{0.3cm}

%%%%%%%%%%%%%%%%%%%%%%%%%%%%%%%%%%%%%%%%%%%%%%%%%%%%%%%%%%%%%%%%%
%%%%%%%%%%%%%%%%%%%%%%%%% subsection 3

{\bf 2.3. Hypertangent divisors.} Getting back to the coordinates
$z_1,\dots,z_N$, write down
$$
h_{\leqslant i}=h_1+\dots+h_i
$$
for $i=1,\dots,\mathop{\rm deg}X$ and consider the second
hypertangent system
$$
\Lambda^R_2=|s_0h_{\leqslant 2}+s_1h_1|_R,
$$
where $s_0\in{\mathbb C}$ and $s_1$ runs through the space of
linear forms in $z_*$. By the condition (N3) the base set
$\mathop{\rm Bs}\Lambda^R_2$ is irreducible and reduced, and by
the condition (N1) it is of codimension 2 on $R$. Therefore, a
general divisor $D_2\in\Lambda^R_2$ does not contain the prime
divisor $D_R$ as a component, so that we get a well defined
effective cycle
$$
Y_2=(D_2\circ D_R)
$$
of codimension 2 on $R$, satisfying the inequality
$$
\frac{\mathop{\rm mult_o}}{\mathop{\rm deg}}Y_2>
\frac{3}{\mathop{\rm deg}X}.
$$
By the linearity of the equivalent inequality
$$
\mathop{\rm mult}\nolimits_oY_2> \frac{3}{\mathop{\rm
deg}X}\mathop{\rm deg}Y_2
$$
in $Y_2$ we may replace the cycle $Y_2$ by its suitable
irreducible component and assume $Y_2$ to be an irreducible
subvariety of codimension 2.\vspace{0.1cm}

{\bf Lemma 2.4.} {\it The subvariety $Y_2$ is not contained in the
tangent divisor $T_R$.}\vspace{0.1cm}

{\bf Proof.} The base set of the hypertangent system $\Lambda^R_2$
is
$$
S_R=\{h_1|_R=h_2|_R=0\}.
$$
It is irreducible, reduced and therefore
$$
\mathop{\rm deg}S_R=2\mathop{\rm deg}X.
$$
By the condition (N1) the equality
$$
\mathop{\rm mult}\nolimits_oS_R=6
$$
holds. Therefore, $Y_2\neq S_R$. However, a certain polynomial
$$
s_0 h_{\leqslant 2}+s_1h_1
$$
vanishes on $Y_2$, where $s_0\neq 0$, since the divisor
$D_2\in\Lambda^R_2$ is chosen to be general. If we had
$$
h_1|_{Y_2}\equiv 0,
$$
then we would have got $h_{\leqslant 2}|_{Y_2}\equiv 0$. Since
$h_{\leqslant 2}=h_1+h_2$, this would have implied that
$h_2|_{Y_2}\equiv 0$ and $Y_2\subset\mathop{\rm
Bs}\Lambda^R_2=S_R$, which is not true. Q.E.D. for the
lemma.\vspace{0.1cm}

By the lemma that we have just shown, the effective cycle
$$
Y_3=(Y_2\circ T_R)
$$
of codimension 3 on $R$ is well defined. It satisfies the
inequality
$$
\frac{\mathop{\rm mult_o}}{\mathop{\rm deg}}Y_3>
\frac{6}{\mathop{\rm deg}X}.
$$
The cycle $Y_3$ can be assumed to be an irreducible subvariety of
codimension 3 on $R$ for the same reason as $Y_2$.\vspace{0.1cm}

Now applying the technique of hypertangent divisors in the usual
way \cite[Chapter 3]{Pukh13a}, we intersect $Y_3$ with general
hypertangent divisors
$$
D_4\in\Lambda^R_4,\,\,\dots,\,\, D_{N-4}\in\Lambda^R_{N-4},
$$
using the condition (N1), and obtain an irreducible curve
$C\subset R$, satisfying by (\ref{19.04.2019.5}) the inequality
$$
\frac{\mathop{\rm mult_o}}{\mathop{\rm deg}}C>
\frac{6}{\mathop{\rm deg}X}\cdot\frac54\cdot\dots\cdot
\frac{N-2}{N-3}=\frac{6}{\mathop{\rm deg}X}\cdot\frac{N-3}{4}
\geqslant 1,
$$
which is impossible.\vspace{0.1cm}

This proves Proposition 2.1. Q.E.D.\vspace{0.1cm}

{\bf Proof of Proposition 1.1.} It is sufficient to check that the
hypersurface $\Gamma$ satisfies all the assumptions that were made
about the hypersurface $X$. Indeed, $\Gamma$ has at most quadratic
singularities of rank $\geqslant 8$, so that
$$
\mathop{\rm codim\,} (\mathop{\rm Sing\,} X\subset X)\geqslant 7.
$$
That the inequality (\ref{19.04.2019.5}) is true for $\Gamma$, one
checks by elementary computations. The condition (N1) follows from
(R1.1), the condition (N2) from (R1.2), the condition (N3) from
(R1.3). Therefore, we can apply Proposition 2.1. This proves
Proposition 1.1. Q.E.D.

%%%%%%%%%%%%%%%%%%%%%%%%%%%%%%%%%%%%%%%%%%%%%%%%%%%%%%%%%%%%%%%%%
%%%%%%%%%%%%%%%%%%%%%%%%%%%%%%%%%%%%%%%%%%%%%%%%%%%%%%%%%%%%%%%%%
%%%%%%%%%%%%%%%%%%%%   SECTION 3

\section{Reduction to a hyperplane section}

In this section we consider hypersurfaces $X\subset{\mathbb P}^N$
with at most quadratic singularities, the rank of which is bounded
from below, which also satisfy some additional conditions. For a
non-canonical pair $(X,\frac{1}{n} D_X)$, where $D_X\sim nH_X$
does not contain hyperplane sections of the hypersurface $X$, we
construct a special hyperplane section $\Delta$, such that the
pair $(\Delta,\frac{1}{n} D_{\Delta})$, where
$D_{\Delta}=D_X|_{\Delta}$, is again non-canonical and, into the
bargain, somewhat ``better'' than the original pair: the
multiplicity of the divisor $D_{\Delta}$ at some point $o\in
\Delta$ is higher than the multiplicity of the original divisor
$D_X$ at this point.\vspace{0.3cm}

{\bf 3.1. Hypersurfaces with singularities.} Take $N\geqslant 8$
and let $X\subset{\mathbb P}^N$ be a hypersurface, satisfying the
following conditions:\vspace{0.1cm}

(S1) every point $o\in X$ is either non-singular, or a quadratic
singularity of rank $\geqslant 7$,\vspace{0.1cm}

(S2) for every effective divisor $D\sim nH_X$, where
$H_X\in\mathop{\rm Pic}X$ is the class of a hyperplane section and
$n\geqslant 1$, the union $\mathop{\rm CS}(X,\frac{1}{n}D_X)$ os
the centres of all non log canonical singularities of the pair
$(X,\frac{1}{n}D_X)$ is contained in $\mathop{\rm
Sing}X$,\vspace{0.1cm}

(S3) for every effective divisor $Y$ on the section of $X$ by a
linear subspace of codimension 1 or 2 in ${\mathbb P}^N$ and every
point $o\in Y$, singular on $X$, the following inequality holds:
\begin{equation}\label{07.06.2019.1}
\frac{\mathop{\rm mult}\nolimits_o}{\mathop{\rm deg}}Y<
\frac{4}{\mathop{\rm deg}X}.
\end{equation}

The condition (S1), Grothendieck's theorem on parafactoriality
\cite{CL} and the Lefschetz theorem imply that $X$ is a factorial
variety and $\mathop{\rm Cl} X=\mathop{\rm Pic} X={\mathbb Z}
H_X$, since $\mathop{\rm codim\,}(\mathop{\rm Sing} X\subset
X)\geqslant 6$. As every hyperplane section of the hypersurface
$X$ is a hypersurface in ${\mathbb P}^{N-1}$, the singular locus
of which has codimension at least 4, it is also
factorial.\vspace{0.1cm}

Assume, furthermore, that $D_X\sim nH_X$ is an effective divisor,
such that we have $\mathop{\rm
CS}(X,\frac{1}{n}D_X)\neq\emptyset$, and moreover, there is a
point $o\in\mathop{\rm CS}(X,\frac{1}{n}D_X)\subset \mathop{\rm
Sing} X$ (see the condition (S2)), which is a quadratic
singularity of rank $\geqslant 8$. Let
$$
\varphi\colon X^+\to X
$$
be its blow up with the exceptional divisor $E=\varphi^{-1}(o)$,
which by our assumption is a quadric of rank $\geqslant 8$. For
the strict transform $D^+_X\subset X^+$ we can write
$$
D^+_X\sim n(H_X-\alpha E),
$$
where by the condition (S3) we have $\alpha<2$, since $\mathop{\rm
mult}\nolimits_oD_X<4n$.\vspace{0.1cm}

{\bf Remark 3.1.} As we will see below, under our assumptions the
inequality $\alpha >1$ holds. Since for every hyperplane section
$\Delta\ni o$ of the hypersurface $X$ and its strict transform
$\Delta^+\subset X^+$ we have
$$
\Delta^+\sim\Delta-E,
$$
the pair $(X,\Delta)$ is canonical, so that we may assume that the
effective divisor $D_X$ does not contain hyperplane sections of
the hypersurface $X$ as components (if there are such components,
they can be removed with all assumptions being kept). For that
reason, for any hyperplane section $\Delta\ni o$ the effective
cycle $(\Delta\circ D_X)$ of codimension 2 on $X$ is well defined.
We will understand this cycle as an effective divisor on the
hypersurface $\Delta\subset{\mathbb P}^{N-1}$ and denote it by the
symbol $D_{\Delta}$.\vspace{0.1cm}

{\bf Proposition 3.1.} {\it There is a hyperplane section
$\Delta\ni o$ of the hypersurface $X$, such that
$$
o\in\mathop{\rm CS}(\Delta,\frac{1}{n}D_{\Delta})
$$
and $D^+_{\Delta}\sim n(H_{\Delta}-\alpha_{\Delta}E_{\Delta})$,
and moreover the following inequality holds:}
\begin{equation}\label{27.03.2019.3}
\alpha_{\Delta}>\frac12\alpha+1
\end{equation}

(Here $H_{\Delta}$ is the class of a hyperplane section of the
hypersurface $\Delta\subset{\mathbb P}^{N-1}$, and
$E_{\Delta}=\Delta^+\cap E$ is the exceptional divisor of the blow
up $\varphi_{\Delta}\colon\Delta^+\to\Delta$, where $\Delta^+$ is
the strict transform of $\Delta$ on $X^+$ and $D^+_{\Delta}$ is
the strict transform of the divisor $D_{\Delta}$ on
$\Delta^+$.)\vspace{0.1cm}

{\bf Proof.} Obviously, $D_{\Delta}\sim nH_{\Delta}$. We have
$$
\mathop{\rm CS}(\Delta,\frac{1}{n}D_{\Delta})\supset
\Delta\cap\mathop{\rm CS}(X,\frac{1}{n}D_X)
$$
for {\it every} hyperplane section $\Delta$, so that we only need
to show the existence of the hyperplane section $\Delta$ for which
the inequality (\ref{27.03.2019.3}) is satisfied. This fact is
obtained by the arguments, repeating the proof of Theorem 1.4 in
\cite{Pukh15a} (Subsections 4.2, 4.3) almost word for word. We
will go through the main steps of these arguments, dwelling on the
necessary modifications. Whenever possible, we use the same
notations as in \cite[Subsections 4.2,
4.3]{Pukh15a}.\vspace{0.3cm}

%%%%%%%%%%%%%%%%%%%%%%%%%%%%%%%%%%%%%%%%%%%%%%%%%%%%%%%%%%%%%%%%
%%%%%%%%%%%%%%%%%%%%%%%%% subsection 2

{\bf 3.2. Preliminary constructions.} Consider the section $P$ of
the hypersurface $X$ by a general 5-dimensional linear subspace,
containing the point $o$. Obviously, $P\subset{\mathbb P}^5$ is a
factorial hypersurface, $o\in P$  is an isolated quadratic
singularity of the maximal rank. Let $P^+\subset X^+$ be the
strict transform of the hypersurface $P$, so that $E_P=P^+\cap E$
is a non-singular three-dimensional quadric. Set $D_P=(D\circ
P)=D|_P$. Obviously, by the inversion of adjunction the pair
$(P,\frac{1}{n}D_P)$ has the point $o$ as an isolated centre of a
non log canonical singularity. Since $a(E_P)=2$ and $D^+_P\sim
nH_P-\alpha nE_P$ (where $H_P$ is the class of a hyperplane
section of the hypersurface $P\subset{\mathbb P^5}$), and moreover
$\alpha<2$, we conclude that the pair $(P^+,\frac{1}{n}D^+_P)$ is
not log canonical and the union
$$
\mathop{\rm LCS}\nolimits_E(P^+,\frac{1}{n}D^+_P)
$$
of the centres of all non log canonical singularities of that
pair, intersecting the exceptional divisor $E_P$, is a connected
closed subset of the quadric $E_P$. Let $S_P$ be an irreducible
component of maximal dimension of that set. Since $S_P$ is the
centre of certain non log canonical singularity of the pair
$(P^+,\frac{1}{n}D^+_P)$, the inequality
$$
\mathop{\rm mult}\nolimits_{S_P}D^+_P>n
$$
holds. Furthermore, $\mathop{\rm codim}(S_P\subset
E_P)\in\{1,2,3\}$ (and if $S_P$ is a point, then we have
$\mathop{\rm LCS}_E(P^+,\frac{1}{n}D^+_P)=S_P$ by the
connectedness of that set). Coming back to the original pair
$(X,\frac{1}{n}D_X)$, we see that the pair
$(X^+,\frac{1}{n}D^+_X)$ has a non log canonical singularity, the
centre of which is an irreducible subvariety $S\subset E$, such
that $S\cap E_P=S_P$; in particular,
$$
\mathop{\rm codim}(S\subset E)= \mathop{\rm codim}(S_P\subset
E_P)\in\{1,2,3\},
$$
and if the last codimension is equal to 3, then $S\cap E_P$ is a
point and for that reason $S\subset E$ is a linear subspace of
codimension 3. However, on a quadric of rank $\geqslant 8$ there
can be no linear subspaces of codimension 3, so that $\mathop{\rm
codim}(S\subset E)\in\{1,2\}$.\vspace{0.1cm}

{\bf Proposition 3.2.} {\it The case $\mathop{\rm codim}(S\subset
E)=1$ is impossible.}\vspace{0.1cm}

{\bf Proof.} Assume that this case takes place. Then $S\subset E$
is a prime divisor, which is cut out on $E$ by a hypersurface of
degree $d_S\geqslant 1$, that is, $S\sim d_SH_E$, where $H_E$ is
the class of a hyperplane section of the quadric $E$. We have
$(D^+_X\circ E)\sim\alpha nH_E$, so that
$$
2>\alpha\geqslant d_S,
$$
and for that reason $S\sim H_E$ is a hyperplane section of the
quadric $E$. Let $\Delta\in|H|$ be the uniquely determined
hyperplane section of the hypersurface $X$, such that $\Delta\ni
o$ and $(\Delta^+\circ E)=\Delta^+\cap E=S$. For the effective
divisor $D_{\Delta}$ the inequality
$$
\mathop{\rm mult}\nolimits_oD_{\Delta}\geqslant 2\alpha
n+2\mathop{\rm mult}\nolimits_SD^+_X>4n
$$
holds. Taking into account that $\mathop{\rm deg}(\Delta\circ
D_X)=n\mathop{\rm deg}X$, we get a contradiction with the
condition (S3), which by assumption is satisfied for the
hypersurface $X$.\vspace{0.1cm}

Q.E.D. for the proposition.\vspace{0.3cm}

%%%%%%%%%%%%%%%%%%%%%%%%%%%%%%%%%%%%%%%%%%%%%%%%%%%%%%%%%%%%%%%%%%
%%%%%%%%%%%%%%%%%%%%%%%%% subsection 3

{\bf 3.3. The case of codimension 2.} We proved above that
$S\subset E$ is a subvariety of codimension 2. Following
\cite[Section 3]{Pukh18b}, for distinct points $p\neq q$ on the
quadric $E$ we denote by the symbol $[p,q]$ the line joining these
two points, {\it provided that it is contained in} $E$, and the
empty set, otherwise, and set
$$
\mathop{\rm Sec\,} (S\subset
E)=\overline{\bigcup_{\begin{array}{l}
p,q\in S\\
p\neq q\\
\end{array}} [p,q]}
$$
(where the line above means the closure).\vspace{0.1cm}

{\bf Lemma 3.1.} {\it One of the following two options takes
place:\vspace{0.1cm}

{\rm (1)} $\mathop{\rm Sec\,}(S\subset E)$ is a hyperplane section
of the quadric $E$, on which $S$ is cut out by a hypersurface of
degree $d_S\geqslant 2$,\vspace{0.1cm}

{\rm (2)} $S=\mathop{\rm Sec\,}(S\subset E)$ is the section of the
quadric $E$ by a linear subspace of codimension 2.}\vspace{0.1cm}

{\bf Proof} repeats the proof of Lemma 4.1 in \cite{Pukh15a}, and
we do not give it here. (The key point in the arguments is that
due to the inequality $\alpha<2$ every line $L=[p,q]\subset E$,
joining some point $p,q\in S$ and lying on $E$, is contained in
$D^+_X$, because $\mathop{\rm mult}_SD^+_X>n$.)\vspace{0.1cm}

{\bf Proposition 3.3.} {\it The option (2) does not take
place.}\vspace{0.1cm}

{\bf Proof.} Assume the converse: the case (2) takes place. Let
$P\subset X$ be the section of the hypersurface $X$ by the linear
subspace of codimension 2 in ${\mathbb P}^N$, that is uniquely
determined by the conditions $P\ni o$ and $P^+\cap
E=S$.\vspace{0.1cm}

The symbol $|H-P|$ stands for the pencil of hyperplane sections of
the hypersurface $X$, containing $P$. For a general divisor
$\Delta\in|H-P|$ we have the equality
$$
\mathop{\rm mult}\nolimits_SD^+_{\Delta}= \mathop{\rm
mult}\nolimits_SD^+_X.
$$

Write down $(\Delta\circ D_X)=G+aP$, where $a\in{\mathbb Z}_+$ and
$G$ is an effective divisor on $\Delta$, not containing $P$ as a
component. Obviously, $G\in|mH_{\Delta}|$, where $m=n-a$ and
$H_{\Delta}$ is the class of a hyperplane section of
$\Delta\subset{\mathbb P}^{N-1}$. The symbols $G^+$ and $\Delta^+$
stand for the strict transforms of $G$ and $\Delta$ on $X^+$,
respectively. Now
$$
G^+\sim mH_{\Delta}-(\alpha n-a)E_{\Delta},
$$
where $E_{\Delta}=\Delta^+\cap E$ is a hyperplane section of the
quadric $E$ and, besides,
$$
\mathop{\rm mult}\nolimits_SG^+= \mathop{\rm
mult}\nolimits_SD^+_X-a>m.
$$

By construction, the effective cycle $(G\circ P)$ of codimension 2
on $\Delta$ is well defined. One can consider it as an effective
divisor on the hypersurface $P\subset{\mathbb P}^{N-2}$. The
following inequality holds:
$$
\mathop{\rm mult}\nolimits_o(G\circ P)\geqslant 2(\alpha n-a)+
2\mathop{\rm mult}\nolimits_SG^+>4m.
$$
Since $\mathop{\rm deg}(G\circ P)=m\mathop{\rm deg}X$, we obtain a
contradiction with the condition (S3), which is satisfied for the
hypersurface $X$. Q.E.D. for the proposition.\vspace{0.3cm}

%%%%%%%%%%%%%%%%%%%%%%%%%%%%%%%%%%%%%%%%%%%%%%%%%%%%%%%%%%%%%%%%%
%%%%%%%%%%%%%%%%%%%%%%%%% subsection 4

{\bf 3.4. The hyperplane section $\Delta$.} We have shown that the
case (1) takes place. Set $\Lambda=S=\mathop{\rm Sec\,}(S\subset
E)$. This is a hyperplane section of the quadric $E$, where
$\Lambda\subset D^+_X$. Set
$$
\mu=\mathop{\rm mult}\nolimits_SD^+_X,\quad\gamma= \mathop{\rm
mult}\nolimits_{\Lambda}D^+_X.
$$
We know that $\mu>n$ and $\mu\leqslant\alpha n<2n$ (the second
inequality holds, because for a general linear subspace
$\Pi\subset E$ of maximal dimension the divisor
$(D^+_X\circ\Pi)=D^+_X\cap\Pi$ on $\Pi$ is a hypersur\-face of
degree $\alpha n$, containing every point of the set $S\cap\Pi$
with multiplicity $\geqslant\mu$).\vspace{0.1cm}

{\bf Lemma 3.2.} {\it The following inequality holds:}
$$
\gamma\geqslant\frac13(2\mu-\alpha n).
$$

{\bf Proof:} this is Lemma 4.2 in \cite{Pukh15a}. The claim of the
lemma is a local fact and for that reason the proof given in
\cite[п. 4.3]{Pukh15a} does not require any modifications and
works word for word. Q.E.D. for the lemma.\vspace{0.1cm}

Now let us consider the uniquely determined hyperplane section
$\Delta$ of the hypersur\-face $X\subset{\mathbb P}^N$, such that
$\Delta\ni o$ and $\Delta^+\cap E=\Lambda$, where $\Delta^+\subset
X^+$ is the strict transform of $\Delta$. Write down
$$
(D^+_X\circ\Delta^+)=D^+_{\Delta}+a\Lambda.
$$
Obviously, $\mathop{\rm mult}_o D_{\Delta}=2(\alpha n+a)$, so that
$$
\alpha_{\Delta}=\alpha+\frac{a}{n}
$$
(recall that $D^+_{\Delta}\sim
n(H_{\Delta}-\alpha_{\Delta}E_{\Delta})$, where
$E_{\Delta}=\Lambda$). Since the subvariety $S$ is cut out on the
quadric $\Lambda$ by a hypersurface of degree $d_S\geqslant 2$, we
obtain the inequality
\begin{equation}\label{15.04.2019.3}
\mathop{\rm mult}\nolimits_SD^+_{\Delta}\leqslant
\frac{1}{d_S}\alpha_{\Delta}n\leqslant\frac{\alpha n+a}{2}.
\end{equation}
Since $S\subset\mathop{\rm LCS}(X^+,\frac{1}{n}D^+_X)$, we get:
$$
S\subset\mathop{\rm
LCS}(\Delta^+,\frac{1}{n}(D^+_{\Delta}+a\Lambda)).
$$
Consider the blow up $\sigma_S\colon\widetilde{\Delta}\to\Delta^+$
of the subvariety $S\subset\Delta^+$ of codimension 2 and denote
its exceptional divisor $\sigma_S^{-1}(S)$ by the symbol
$E_S$.\vspace{0.1cm}

{\bf Proposition 3.4.} {\it For some irreducible divisor
$S_1\subset E_S$, such that the projection $\sigma_S|_{S_1}$ is
birational, the inequality
\begin{equation}\label{15.04.2019.2}
\mathop{\rm mult}\nolimits_S(D^+_{\Delta}+a\Lambda)+ \mathop{\rm
mult}\nolimits_{S_1}(\widetilde{D}_{\Delta}+
a\widetilde{\Lambda})>2n
\end{equation}
holds, where $\widetilde{D}_{\Delta}$ and $\widetilde{\Lambda}$
are the strict transforms, respectively, of $D^+_{\Delta}$ and
$\Lambda$ on $\widetilde{\Delta}$.}\vspace{0.1cm}

{\bf Proof.} This is a well known fact, see \cite[Proposition
9]{Pukh05}. (Note that the subvariety $S$ is, generally speaking,
singular, however $\Delta^+$ is non-singular at the general point
of $S$ and $\widetilde{\Delta}$ is non-singular at the general
point of $S_1$.)\vspace{0.3cm}

%%%%%%%%%%%%%%%%%%%%%%%%%%%%%%%%%%%%%%%%%%%%%%%%%%%%%%%%%%%%%%%%%%%
%%%%%%%%%%%%%%%%%%%%%%%%% subsection 5

{\bf 3.5. End of the proof.} Set $\mu_S=\mathop{\rm
mult}_SD^+_{\Delta}$ and $\beta=\mathop{\rm
mult}_{S_1}\widetilde{D}_{\Delta}$. One of the two cases takes
place:\vspace{0.1cm}

--- the case of general position $S_1\neq E_S\cap\widetilde{\Lambda}$,
so that $S_1\not\subset\widetilde{\Lambda}$,\vspace{0.1cm}

--- the special case
$S_1=E_S\cap\widetilde{\Lambda}$.\vspace{0.1cm}

Let us consider them separately. In the case of general position
the inequality (\ref{15.04.2019.2}) takes the form
$$
\mu_S+\beta+a>2n,
$$
since $\mathop{\rm mult}_{S_1}\widetilde{\Lambda}=0$. Furthermore,
$\mu_S\geqslant\beta$, so that the more so
$$
2\mu_S+a>2n.
$$
On the other hand, from the inequality (\ref{15.04.2019.3}) we get
$2\mu_S\leqslant\alpha n+a$, which implies that
$$
\alpha n+2a>2n
$$
and for that reason
$$
2\alpha_{\Delta}n=2\alpha n+2a>(\alpha+2)n.
$$
The inequality (\ref{27.03.2019.3}) in the case of general
position is now proven.\vspace{0.1cm}

Let us consider the special case. Here $\mathop{\rm
mult}_{S_1}\widetilde{\Lambda}=1$, so that the inequality
(\ref{15.04.2019.2}) takes the form
$$
\mu_S+\beta+2a>2n.
$$
Besides, the effective cycle $(D^+_{\Delta}\circ\Lambda)$,
considered as an effective divisor on $\Lambda$, is cut out on the
quadric $\Lambda$ by a hypersurface of degree $\alpha n+a$, and
contains the divisor $S\sim d_SH_{\Lambda}$ (where $H_{\Lambda}$
is the class of a hyperplane section of $\Lambda$) with
multiplicity $\geqslant\mu_S+\beta$, so that
$$
2(\mu_S+\beta)\leqslant\alpha n+a,
$$
whence we get $\alpha n+5a>4n$ and for that reason
$$
5\alpha_{\Delta} n=5(\alpha n+a)>4(\alpha+1)n,
$$
that is, $\alpha_{\Delta}>\frac45\alpha+\frac45>\frac35\alpha+1$
(since $\alpha>1$). This inequality is stronger than
(\ref{27.03.2019.3}), which completes the proof in the special
case.\vspace{0.1cm}

Q.E.D. for Proposition 3.1.\vspace{0.1cm}

{\bf Proof of Proposition 1.2.} Let us check that the operation of
reduction, described in Subsection 3.1, can be $\rho$ times
applied to the hypersurface $\Gamma\subset {\mathbb P}^M$.
Consider the hypersurface $\Gamma_i\subset {\mathbb P}^{M-i}$,
where $i\in \{0,\dots,\rho-1\}$. Let us show, in the first place,
that $\Gamma_i$ satisfies the condition (S1). Let $p\in \Gamma_i$
be an arbitrary singularity. If $i=0$, then by the condition
(R2.1), the point $p$ is a quadratic singularity of rank
$\geqslant 8$. If $i\geqslant 1$, then there are two options:
either $p\in\Gamma$ is a non-singular point, or $p\in\Gamma$ is a
singularity (recall that $\Gamma_i$ is a section of the
hypersurface $\Gamma$ by a linear subspace of codimension $i$ in
${\mathbb P}^M$). In the second case by the condition (R2.1) the
point $p$ is a quadratic singularity of $\Gamma$ of rank
$\geqslant 2\rho+6\geqslant 2i+8$, since $\rho\geqslant i+1$.
Since a hyperplane section of a quadric of rank $r\geqslant 3$ is
a quadric of rank $\geqslant r-2$, we conclude that $p\in
\Gamma_i$ is a quadratic singularity of rank $\geqslant 8$, so
that the condition (S1) is satisfied at that point (for the
hypersurface $\Gamma_i$).\vspace{0.1cm}

In the first case the point $p$ is non-singular on $\Gamma$, so
that $\Gamma_i$ is a section of $\Gamma$ by a linear subspace of
codimension $i$, which is contained in the tangent hyperplane
$T_p\Gamma$. By the condition (R1.4) the point $p\in\Gamma_i$ is a
quadratic singularity of rank
$$
\geqslant 8+2(i+1)-4-2(i-1)=8
$$
(one should take into account that the cutting subspace is of
codimension $i-1$ in $T_p\Gamma$). Therefore, the condition (S1)
is satisfied in any case.\vspace{0.1cm}

Let us show that the hypersurface $\Gamma_i$ satisfies the
condition (S2) as well. In order to do it, we must check that for
$\Gamma_i$ all assumptions of Subsection 2.1 are satisfied. By
what was said above, the codimension of the set $\mathop{\rm
Sing\,}\Gamma_i$ with respect to $\Gamma_i$ it at least 7
--- this is higher than we need. The inequality (\ref{19.04.2019.5})
takes the form of the estimate
$$
M-i-2\leqslant dl\leqslant \frac32(M-i-3),
$$
which is easy to check. Finally, the conditions (N1), (N2) and
(N3) follow from the conditions (R1.1), (R1.2) and (R1.3),
respectively. By Proposition 2.1 we conclude that the hypersurface
$\Gamma_i$ satisfies the condition (S2).\vspace{0.1cm}

Finally, let us consider the condition (S3). Obviously, it is
sufficient to check that the inequality (\ref{07.06.2019.1}) holds
for any prime divisor $Y$ on the section of the hypersurface
$\Gamma_i$ by a linear subspace $P^*$ of codimension 2 in
${\mathbb P}^{M-i}$. Assume the converse:
\begin{equation}\label{07.06.2019.2}
\frac{\mathop{\rm mult}\nolimits_o}{\mathop{\rm deg}}Y<
\frac{4}{dl}.
\end{equation}
In some affine coordinates with the origin at the point $o$ on the
subspace $P^*={\mathbb P}^{M-i-2}$ the equation of the
hypersurface $P^*\cap\Gamma_i$ has the form
$$
0=q^*_2+q^*_3+\cdots +q^*_{dl},
$$
where by the condition (R2.2) the sequence of homogeneous
polynomials
$$
q^*_2,\quad q^*_3,\quad\dots,\quad q^*_{M-\rho-4}
$$
is regular. Consider general hypertangent divisors
$$
D^*_2,\quad D^*_3,\quad\dots,\quad D^*_{M-\rho-5}.
$$
The first hypertangent divisor $D^*_2$ is irreducible and
satisfies the equality
$$
\frac{\mathop{\rm mult}\nolimits_o}{\mathop{\rm deg}}
D^*_2=\frac{3}{dl},
$$
so that $D^*_2\neq Y$ and the effective cycle of the
scheme-theoretic intersection $(D^*_2\circ Y)$ is well defined. It
satisfies the inequality
$$
\frac{\mathop{\rm mult}\nolimits_o}{\mathop{\rm deg}} (D^*_2\circ
Y) >\frac{6}{dl}.
$$
Let $Y_2$ be an irreducible component of that cycle with the
maximal value of the ratio $\mathop{\rm
mult}\nolimits_o/\mathop{\rm deg}$. Intersecting $Y_2$ with the
divisors
$$
D^*_4,\quad\dots, \quad D^*_{M-\rho-5}
$$
in the usual way (see \cite[Chapter 3]{Pukh13a} or Subsection 2.3
of the present paper), we construct a sequence of irreducible
subvarieties
$$
Y^*_3,\quad\dots,\quad Y^*_{M-\rho-6},
$$
where $\mathop{\rm codim\,}(Y^*_j\subset (P^*\cap\Gamma_i))=j$ and
the last subvariety $Y^*_{M-\rho-6}$ (the dimension of which is
$\rho-i+3\geqslant 4$) satisfies the inequality
$$
\frac{\mathop{\rm mult}\nolimits_o}{\mathop{\rm deg}}
Y^*_{M-\rho-6}
>\frac{6}{dl}\cdot\frac54\cdot\,\cdots\,\cdot\frac{M-\rho-4}{M-\rho-5}=
\frac{3(M-\rho-4)}{2dl}.
$$
It is easy to check that the right hand side of the last
inequality for the values of $d$, $l$ and $\rho$ under
consideration is higher than 1, which gives a contradiction with
the assumption (\ref{07.06.2019.2}) and proves that the
hypersurface $\Gamma_i$ satisfies the condition
(S3).\vspace{0.1cm}

Note that all singular points of $\Gamma_i$ are quadratic
singularities of rank $\geqslant 8$, so that the additional
assumption about the point $o$ made in Subsection 3.1 is
satisfied.\vspace{0.1cm}

Now applying Proposition 3.1, we complete the proof of Proposition
1.2.

\begin{flushleft}
Department of Mathematical Sciences,\\
The University of Liverpool
\end{flushleft}

\noindent{\it pukh@liverpool.ac.uk}

\end{document}